\documentclass[11pt]{article}
\usepackage{amssymb}
\usepackage{amsmath}
\usepackage[all]{xy}
\usepackage{times}
\usepackage{graphicx}

\title{A separation result for countable unions of Borel rectangles}
\author{Dominique LECOMTE}
\date{\today}

\def\ufootnote#1{\let\savedthfn\thefootnote\let\thefootnote\relax
\footnote{#1}\let\thefootnote\savedthfn\addtocounter{footnote}{-1}}

\newcommand{\Ana}{{\it\Sigma}^{1}_{1}}

\newcommand{\Boraone}{{\it\Sigma}^{0}_{1}}

\newcommand{\Borel}{{\it\Delta}^{1}_{1}}

\newcommand{\boraone}{{\bf\Sigma}^{0}_{1}}
\newcommand{\boratwo}{{\bf\Sigma}^{0}_{2}}

\newcommand{\boraxi}{{\bf\Sigma}^{0}_{\xi}}

\newcommand{\bormone}{{\bf\Pi}^{0}_{1}}

\newcommand{\bormtwo}{{\bf\Pi}^{0}_{2}}

\newcommand{\bormlxi}{{\bf\Pi}^{0}_{<\xi}}

\newcommand{\bormxi}{{\bf\Pi}^{0}_{\xi}}

\newcommand{\borxi}{{\bf\Delta}^{0}_{\xi}}
\newcommand{\borme}{{\bf\Pi}^{0}_{\eta}}
\newcommand{\borae}{{\bf\Sigma}^{0}_{\eta}}

\newcommand{\boraep}{{\bf\Sigma}^{0}_{\eta +1}}

\newtheorem{thm} {Theorem} [section]

\newtheorem{lem} [thm] {Lemma}

\newtheorem{them} {Theorem} [subsection]
\newtheorem{defin} [them] {Definition}

\newtheorem{lemm} [them] {Lemma}

\begin{document}

\maketitle

\centerline{$\bullet$ Universit\' e Paris 6, Institut de Math\'ematiques de Jussieu-Paris Rive Gauche,}

\centerline{Projet Analyse Fonctionnelle, Couloir 16-26, 4\`eme \'etage, Case 247,}

\centerline{4, place Jussieu, 75 252 Paris Cedex 05, France}

\centerline{dominique.lecomte@upmc.fr}\medskip

\centerline{$\bullet$ Universit\'e de Picardie, I.U.T. de l'Oise, site de Creil,}

\centerline{13, all\'ee de la fa\"\i encerie, 60 107 Creil, France}\medskip\medskip\medskip\medskip\medskip\medskip

\ufootnote{{\it 2010 Mathematics Subject Classification.}~Primary: 03E15, Secondary: 54H05}

\ufootnote{{\it Keywords and phrases.}~Borel class, countable Borel coloring, Borel rectangle}

\noindent {\bf Abstract.} We provide dichotomy results characterizing when two disjoint analytic binary relations can be separated by a countable union of $\boraone\!\times\!\boraxi$ sets, or by a 
$\bormone\!\times\!\bormxi$ set.

\vfill\eject

\section{$\!\!\!\!\!\!$ Introduction}\indent

 The reader should see [K] for the standard descriptive set theoretic notation and material used in this paper. All our relations will be binary. The motivation for this work goes back to the following so called $\mathbb{G}_0$-dichotomy, essentially proved in [K-S-T]. 
  
\begin{thm} (Kechris, Solecki, Todor\v cevi\' c) \label{G0} There is a Borel relation $\mathbb{G}_0$ on $2^\omega$ such that, for any Polish space $X$ and any analytic relation $A$ on $X$, exactly one of the following holds:\smallskip

(a) there is $c\! :\! X\!\rightarrow\!\omega$ Borel such that $c(x)\!\not=\! c(y)$ if $(x,y)\!\in\! A$ (a {\bf countable Borel coloring} of $A$),\smallskip

(b) there is $f\! :\! 2^\omega\!\rightarrow\! X$ continuous such that $\mathbb{G}_0\!\subseteq\! (f\!\times\! f)^{-1}(A)$.\end{thm}

 This result had a lot of developments since. For instance, Miller developed some techniques to recover many dichotomy results of descriptive set theory, without using effective descriptive set theory (see [M]). He replaces it with some versions of Theorem \ref{G0}. In particular, he can prove Theorem \ref{G0} without effective descriptive set theory. In [L1], the author derives from Theorem \ref{G0} a dichotomy result characterizing when two disjoint analytic sets can be separated by a countable union of Borel rectangles. In order to state it, we give some notation that will also be useful to state our main results.\medskip 
 
\noindent\bf Notation.\rm\ Let, for $\varepsilon\!\in\! 2\! :=\!\{ 0,1\}$, $X_\varepsilon ,Y_\varepsilon$ be Polish spaces, and $A_\varepsilon,B_\varepsilon$ be disjoint analytic subsets of $X_\varepsilon\!\times\! Y_\varepsilon$. We set\medskip

\leftline{$(X_0,Y_0,A_0,B_0)\leq (X_1,Y_1,A_1,B_1)\Leftrightarrow$}\smallskip

\rightline{$\exists f\! :\! X_0\!\rightarrow\! X_1~~\exists g\! :\! Y_0\!\rightarrow\! Y_1\mbox{ continuous with }
A_0\!\subseteq\! (f\!\times\! g)^{-1}(A_1)$ and $B_0\!\subseteq\! (f\!\times\! g)^{-1}(B_1).$}\medskip

\noindent If $X$ is a set, then the {\bf diagonal} of $X$ is $\Delta (X)\! :=\!\{ (x,x)\mid x\!\in\! X\}$. 

\begin{thm} \label{potop} Let $X,Y$ be Polish spaces, and $A,B$ be disjoint analytic subsets of $X\!\times\! Y$. Exactly one of the following holds:\smallskip  

(a) the set $A$ can be separated from $B$ by a countable union of Borel rectangles,\smallskip  

(b) $\big( 2^\omega ,2^\omega ,\Delta (2^\omega ),\mathbb{G}_0\big)\leq (X,Y,A,B)$.\end{thm}

 It is easy to check that Theorem \ref{G0} is also an easy consequence of Theorem \ref{potop}. This means that the study of the countable Borel colorings is highly related to the study of countable unions of Borel rectangles. It is natural to ask for level by level versions of these two results, with respect to the Borel hierarchy. This work was initiated in [L-Z], where the authors prove the following.

\begin{thm} \label{LZ1} (Lecomte, Zelen\'y) Let $\xi\!\in\!\{ 1,2,3\}$. Then we can find a zero-dimensional Polish space $\mathbb{X}$, and an analytic relation $\mathbb{A}$ on $\mathbb{X}$ such that for any (zero-dimensional if $\xi\! =\! 1$) Polish space $X$, and for any relation $A$ on $X$, exactly one of the following holds:\smallskip  

(a) there is a countable $\borxi$-measurable coloring of $A$,\smallskip  

(b) there is $f\! :\!\mathbb{X}\!\rightarrow\! X$ continuous such that $\mathbb{A}\!\subseteq\! (f\!\times\! f)^{-1}(A)$.\end{thm}

 In [L-Z], the authors note that the study of countable $\borxi$-measurable colorings is highly related to the study of countable unions of $\boraxi$ rectangles, since the existence of a countable $\borxi$-measurable coloring of a relation $A$ on a (zero-dimensional if $\xi\! =\! 1$) Polish space $X$ is equivalent to the fact that $\Delta (X)$ can be separated from $A$ by a countable union of $\boraxi$ rectangles, by the generalized reduction property for the class $\boraxi$ (see 22.16 in [K]). In this direction, they prove the following.
 
\vfill\eject
 
\begin{thm} \label{LZ2} (Lecomte, Zelen\'y) Let $\xi\!\in\!\{ 1,2\}$. Then we can find zero-dimensional Polish spaces 
$\mathbb{X},\mathbb{Y}$, and disjoint analytic subsets $\mathbb{A},\mathbb{B}$ of 
$\mathbb{X}\!\times\!\mathbb{Y}$ such that for any Polish spaces $X,Y$, and for any pair $A,B$ of disjoint analytic subsets of $X\!\times\! Y$, exactly one of the following holds:\smallskip  

(a) the set $A$ can be separated from $B$ by a $(\boraxi\!\times\!\boraxi )_\sigma$ set,\smallskip  

(b) $(\mathbb{X},\mathbb{Y},\mathbb{A},\mathbb{B})\leq (X,Y,A,B)$.\end{thm}

 In fact, we can think of a number of related problems of this kind. We can study\medskip
 
\noindent - the finite or bounded finite Borel colorings,\smallskip

\noindent - the separation of disjoint analytic sets by a finite or bounded finite union of Borel rectangles,\smallskip

\noindent - the finite, bounded finite, or infinite Borel colorings of bounded complexity,\smallskip

\noindent - the separation of disjoint analytic sets by a finite, bounded finite or infinite union of Borel rectangles of bounded complexity...\medskip

 This last question has been studied in [Za] in the case of one rectangle. In [Za], the author characterizes when two disjoint analytic sets can be separated by a $\boraone$ (or $\bormxi$ when 
$\xi\!\leq\! 2$) rectangle. Louveau suggested that it could be very interesting to study the non-symmetric version of the problem to understand it better (we can also make this remark for countable unions of rectangles, which is another motivation for Theorem \ref{countunionrectsigma} to come). Zamora noticed that the problems of the separation of analytic sets by a $\bormone\!\times\!\bormtwo$ set and by a $(\boraone\!\times\!\boratwo )_\sigma$ set are very much related (he derives a dichotomy for the rectangles from a dichotomy for the countable unions of rectangles). His technique cannot be extended to higher levels since it uses countability. However, the relation just mentioned is much stronger than in [Za], as we will see. The main results in this paper generalize these two Zamora results, and are, hopefully, steps towards the generalization of Theorem \ref{LZ2}, and then Theorem \ref{LZ1}. The first one is about countable unions of rectangles of the form $\boraone\!\times\!\boraxi$.

\begin{thm} \label{countunionrectsigma} Let $\xi\!\geq\! 1$ be a countable ordinal. Then there are zero-dimensional Polish spaces $\mathbb{X},\mathbb{Y}$, and disjoint analytic subsets 
$\mathbb{A},\mathbb{B}$ of $\mathbb{X}\!\times\!\mathbb{Y}$ such that for any Polish spaces $X,Y$, and for any pair $A,B$ of disjoint analytic subsets of $X\!\times\! Y$, exactly one of the following holds:\smallskip  

(a) the set $A$ can be separated from $B$ by a $(\boraone\!\times\!\boraxi )_\sigma$ set,\smallskip  

(b) $(\mathbb{X},\mathbb{Y},\mathbb{A},\mathbb{B})\leq (X,Y,A,B)$.\end{thm}

The second one is about rectangles of the form $\bormone\!\times\!\bormxi$.

\begin{thm} \label{rectpi} Let $\xi\!\geq\! 1$ be a countable ordinal. Then there are zero-dimensional Polish spaces $\mathbb{X},\mathbb{Y}$, and disjoint analytic subsets $\mathbb{A},\mathbb{B}$ of 
$\mathbb{X}\!\times\!\mathbb{Y}$ such that for any Polish spaces $X,Y$, and for any pair $A,B$ of disjoint analytic subsets of $X\!\times\! Y$, exactly one of the following holds:\smallskip  

(a) the set $A$ can be separated from $B$ by a $\bormone\!\times\!\bormxi$ set,\smallskip  

(b) $(\mathbb{X},\mathbb{Y},\mathbb{A},\mathbb{B})\leq (X,Y,A,B)$.\end{thm}

 One of our key tools to prove these two results is the representation theorem for Borel sets by Debs and Saint Raymond. A classical result of Lusin-Souslin asserts that any Borel subset $\cal B$ of a Polish space is the bijective continuous image of a closed subset of the Baire space (see 13.7 in [K]). There is a level by level version of this result due to Kuratowski: the Baire class of the inverse map of the bijection is essentially equal to the Borel rank of $\cal B$ (see Theorem 1 in [Ku]).
 
\vfill\eject
 
  The representation theorem for Borel sets by Debs and Saint Raymond refines this Kuratowski result (see Theorem I-6.6 in [D-SR]). We will state it and recall the material needed to state it in the next section. Initially, the representation theorem had three applications in [D-SR]: a theorem about continuous liftings, another one about compact covering maps, and a new proof (involving games as in the original paper) of the Louveau-Saint Raymond dichotomy characterizing when two disjoint analytic sets can be separated by a $\boraxi$ (or $\bormxi$) set (see page 433 in [Lo-SR]). In [L3] and [L4], the representation theorem is used to prove a dichotomy about potential Wadge classes. Its proof provides another new proof of the Louveau-Saint Raymond theorem which does not involve games.\medskip

 A very remarkable phenomenon happens in the present paper. In the applications just mentioned, the representation theorem was used only inside the proofs. Here, the representation theorem is used not only in the proofs of Theorems \ref{countunionrectsigma} and \ref{rectpi}, but also to define the minimal objects 
$\mathbb{X},\mathbb{Y},\mathbb{A},\mathbb{B}$. We believe that the minimal objects cannot be that simple for higher levels. Moreover, Theorem \ref{LZ2} provides an extension of Theorem \ref{countunionrectsigma} to countable unions of $\boratwo$ rectangles. It is possible to prove such an  extension using the representation theorem. However, we could not prove further extensions, leaving the general case of countable unions of rectangles of the form $\borae\!\times\!\boraxi$, or just $\boraxi\!\times\!\boraxi$,  open for future work.\medskip
 
 The organization of the paper is as follows. In Section 2, we recall the material about representation needed here, as well as some lemmas from [L3], and we give some effective facts needed to prove our main results. We prove Theorem \ref{countunionrectsigma} in Section 3, and Theorem \ref{rectpi} in Section 4.

\section{$\!\!\!\!\!\!$ Preliminaries}

\subsection{$\!\!\!\!\!\!$ Representation of Borel sets}\indent

 The following definition can be found in [D-SR].

\begin{defin} (Debs-Saint Raymond) A partial order relation $R$ on $2^{<\omega}$ is a 
{\bf tree relation} if, for $s\!\in\! 2^{<\omega}$,\smallskip

(a) $\emptyset ~R~s$,\smallskip

(b) the set $P_{R}(s)\! :=\!\{ t\!\in\! 2^{<\omega}\mid t~R~s\}$ is finite and linearly ordered by $R$ ($h_R(s)$ will denote the number of strict $R$-predecessors of $s$, so that $h_R(s)\! =\!\mbox{Card}\big( P_R(s)\big)\! -\! 1$).\smallskip

\noindent $\bullet$ Let $R$ be a tree relation. An $R$-{\bf branch} is a $\subseteq$-maximal subset of $2^{<\omega}$ linearly ordered by $R$. We denote by $[R]$ the set of all {\rm infinite} $R$-branches.\smallskip

 We equip $(2^{<\omega})^\omega$ with the product of the discrete topology on $2^{<\omega}$. If $R$ is a tree relation, then the space 
$[R]\!\subseteq\! (2^{<\omega})^\omega$ is equipped with the topology induced by that of $(2^{<\omega})^\omega$, and is a Polish space. A basic clopen set is of the form $N^R_s\! :=\!\big\{\gamma\!\in\! [R]\mid\gamma\big(h_R(s)\big)\! =\! s\big\}$, where $s\!\in\! 2^{<\omega}$.\smallskip

\noindent $\bullet$ Let $R$, $S$ be tree relations with $R\!\subseteq\! S$. The 
{\bf canonical map} $\Pi\! :\! [R]\!\rightarrow\! [S]$ is defined by
$$\Pi (\gamma )\! :=\!\mbox{ the unique }S\mbox{-branch containing }\gamma .$$
The canonical map is continuous.

\vfill\eject

\noindent $\bullet$ Let $S$ be a tree relation. We say that $R\!\subseteq\! S$ is {\bf distinguished} in $S$ if
$$\forall s,t,u\!\in\! 2^{<\omega}\ \ \left.
\begin{array}{ll}
& s~S~t~S~u\cr & \cr
& \ \ s~R~u
\end{array}\right\} ~\Rightarrow ~s~R~t.$$
$\bullet$ Let $\eta\! <\!\omega_1$. A family $(R^\rho )_{\rho\leq\eta}$ of tree relations is a 
{\bf resolution family} if\smallskip

(a) $R^{\rho +1}$ is a distinguished subtree of $R^\rho$, for each $\rho\! <\!\eta$.\smallskip

(b) $R^\lambda\! =\!\bigcap_{\rho <\lambda}~R^\rho$, for each limit ordinal $\lambda\!\leq\!\eta$.
\end{defin}

 The representation theorem of Borel sets is as follows in the successor case (see Theorems I-6.6 and I-3.8 in [D-SR]).

\begin{them} \label{rep} (Debs-Saint Raymond) Let $\eta$ be a countable ordinal, and $P\!\in\! {\bf\Pi}^0_{\eta +1}([\subseteq ])$. Then there is a resolution family $(R^\rho )_{\rho\leq\eta}$ such that\smallskip

(a) $R^0\! =\subseteq$,\smallskip

(b) the canonical map $\Pi\! :\! [R^\eta ]\!\rightarrow\! [R^0]$ is a continuous bijection with $\boraep$-measurable inverse,\smallskip

(c) the set $\Pi^{-1}(P)$ is a closed subset of $[R^\eta ]$.\end{them}

 For the limit case, we need some more definition that can be found in [D-SR].

\begin{defin} (Debs-Saint Raymond) Let $\xi$ be an infinite limit countable ordinal. We say that a resolution family $(R^\rho )_{\rho\leq\xi}$ with $R^0\! =\subseteq$ is {\bf uniform} if
$$\forall k\!\in\!\omega ~~\exists\xi_k\! <\!\xi ~~\forall s,t\!\in\! 2^{<\omega}~~
\Big(\mbox{min}\big( h_{R^\xi}(s),h_{R^\xi}(t)\big)\!\leq\! k\wedge s\ R^{\xi_k}\ t\Big)\Rightarrow s\ R^\xi\ t.$$
We may (and will) assume that $\xi_k\!\geq\! 1$.\end{defin}

 The representation theorem of Borel sets is as follows in the limit case (see Theorems I-6.6 and I-4.1 in [D-SR]).
 
\begin{them} \label{replim} (Debs-Saint Raymond) Let $\xi$ be an infinite limit countable ordinal, and $P\!\in\!\bormxi ([\subseteq ])$. Then there is a uniform resolution family 
$(R^\rho )_{\rho\leq\xi}$ such that\smallskip

(a) $R^0\! =\subseteq$,\smallskip

(b) the canonical map $\Pi\! :\! [R^\xi ]\!\rightarrow\! [R^0]$ is a continuous bijection with $\boraxi$-measurable inverse,\smallskip

(c) the set $\Pi^{-1}(P)$ is a closed subset of $[R^\xi ]$.\end{them}
 
 We will use the following extension of the property of distinction (see Lemma 2.3.2 in [L3]).

\begin{lemm} \label{extdist} Let $\eta\! <\!\omega_1$, $(R^\rho )_{\rho\leq\eta}$ be a resolution family, and 
$\rho\! <\!\eta$. Assume that $s,t,u\!\in\! 2^{<\omega}$, $s~R^0~t~R^\rho ~u$ and $s~R^{\rho +1}~u$. Then 
$s~R^{\rho +1}~t$.\end{lemm}

\noindent\bf Notation.\rm ~Let $\eta\! <\!\omega_1$, $(R^\rho )_{\rho\leq\eta}$ be a resolution family with $R^0\! =\subseteq$, $s\!\in\! 2^{<\omega}$, and 
$\rho\!\leq\!\eta$. We define
$$s^\rho\! :=\!\left\{\!\!\!\!\!\!\!\!
\begin{array}{ll}
 & \emptyset\mbox{ if }s\! =\!\emptyset\mbox{,}\cr
 & s\vert\mbox{max}\{ l\! <\!\vert s\vert\mid s\vert l~R^\rho ~s\}\mbox{ if }s\!\not=\!\emptyset .
\end{array}
\right.$$ 
We enumerate $\{ s^\rho\mid\rho\!\leq\!\eta\}$ by $\{ s^{\xi_i}\mid 1\!\leq\! i\!\leq\! n\}$, where 
$n\!\geq\! 1$ is a natural number and $\xi_1\! <\! ...\! <\!\xi_n\! =\!\eta$. We can write 
$s^{\xi_n}\!\subsetneqq\! s^{\xi_{n-1}}\!\subsetneqq\! ...\!\subsetneqq\! s^{\xi_2}\!\subsetneqq\! 
s^{\xi_1}\!\subseteq\! s$. By Lemma \ref{extdist}, $s^{\xi_{i+1}}~R^{\xi_i+1}~s^{\xi_i}$ if 
$1\!\leq\! i\! <\! n$.

\vfill\eject

 We will also use the following lemma (see Lemma 2.3.3 in [L3]).

\begin{lemm} \label{enum} Let $\eta\! <\!\omega_1$, $(R^\rho )_{\rho\leq\eta}$ be a resolution family with $R^0\! =\subseteq$, 
$s\!\in\! 2^{<\omega}\!\setminus\!\{\emptyset\}$ and $1\!\leq\! i\! <\! n$. Then we may assume that $s^{\xi_i+1}\!\subsetneqq\! s^{\xi_i}$.\end{lemm}

\noindent\bf Notation.\rm ~The map $h\! :\! 2^\omega\!\rightarrow\! [\subseteq]$, for which 
$h (\alpha )$ is the strictly $\subseteq$-increasing sequence of all initial segments of $\alpha$, is a homeomorphism.

\subsection{$\!\!\!\!\!\!$ Topologies}\indent

The reader should see [Mo] for the basic notions of effective descriptive set theory.\medskip

\noindent\bf Notation.\rm\ Let $S$ be a recursively presented Polish space.\medskip

\noindent (1) The {\bf Gandy-Harrington topology} on $S$ is generated by $\Ana (S)$ and denoted 
${\it\Sigma}_S$. Recall the following facts about ${\it\Sigma}_S$ (see [L2]).\smallskip

\noindent - ${\it\Sigma}_S$ is finer than the initial topology of $S$.\smallskip

\noindent - We set $\Omega_S :=\{ s\!\in\! S\mid\omega_1^s\! =\!\omega_1^{\mbox{CK}}\}$. Then 
$\Omega_S$ is $\Ana (S)$ and dense in $(S,{\it\Sigma}_S)$.\smallskip

\noindent - $W\cap\Omega_S$ is a clopen subset of $(\Omega_S,{\it\Sigma}_S)$ for each 
$W\!\in\!\Ana (S)$.\smallskip

\noindent - $(\Omega_S,{\it\Sigma}_S)$ is a zero-dimensional Polish space. So we fix a complete compatible metric on $(\Omega_S,{\it\Sigma}_S)$.\medskip

\noindent (2) We call $T_1$ the usual topology on $S$, and $T_\eta$ is the topology generated by the 
$\Ana\cap {\bf\Pi}^0_{<\eta}$ subsets of $S$ if $2\!\leq\!\eta\! <\!\omega_1^{\mbox{CK}}$ (see Definition 1.5 in [Lo]).\medskip

 The next result is essentially Lemma 2.2.2 and the claim in the proof of Theorem 2.4.1 in [L3].

\begin{lemm} \label{top} Let $S$ be a recursively presented Polish space, and 
$1\!\leq\!\eta\! <\!\omega^{\mbox{CK}}_1$.\smallskip

(a) (Louveau) Fix $A\!\in\!\Ana (S)$. Then $\overline{A}^{T_\eta}$ is $\borme$, $\Ana$, and 
$\boraone (T_{\eta +1})$.\smallskip

(b) Let $p\!\geq\! 1$ be a natural number, 
$1\!\leq\!\eta_{1}\! <\!\eta_{2}\! <\!\ldots\! <\!\eta_{p}\!\leq\!\eta$, $S_1$, $\ldots$, $S_p\!\in\!\Ana (S)$, and $O\!\in\!\Boraone (S)$. Assume that $S_i\!\subseteq\!\overline{S_{i+1}}^{T_{\eta_i+1}}$ if 
$1\!\leq\! i\! <\! p$. Then ${S_p\cap\bigcap_{1\leq i<p}~\overline{S_i}^{T_{\eta_i}}}\cap O$ is $T_1$-dense in $\overline{S_1}^{T_1}\cap O$.\smallskip

(c) Let $(R^\rho )_{\rho\leq\eta}$ be a resolution family with $R^0\! =\subseteq$, 
$s\!\in\! 2^{<\omega}\!\setminus\!\{\emptyset\}$, $S_{s^\rho}\!\in\!\Ana (S)$ (for 
$1\!\leq\!\rho\!\leq\!\eta$),  $E\!\in\!\Ana (S)$, and $O\!\in\!\Boraone (S)$. We assume that 
$S_{s^\eta}\!\subseteq\!\overline{E}^{T_{\eta +1}}$ and 
$S_t\!\subseteq\!\overline{S_u}^{T_\rho}\mbox{ if }u~R^\rho ~t\!\subsetneqq\! s$ and 
$1\!\leq\!\rho\!\leq\!\eta$. Then 
$S_{s^\eta}\cap\bigcap_{1\leq\rho <\eta}~\overline{S_{s^\rho}}^{T_\rho}\cap O$ and 
$E\cap\bigcap_{1\leq\rho\leq\eta}~\overline{S_{s^\rho}}^{T_\rho}\cap O$ are $T_1$-dense in 
$\overline{S_{s^1}}^{T_1}\cap O$.\end{lemm}

\noindent\bf Proof.\rm ~(a) See Lemma 1.7 in [Lo].\medskip

\noindent (b) Let $D$ be a $\Boraone$ subset of $S$ meeting $\overline{S_1}^{T_1}\cap O$. Then 
$S_1\cap D\cap O\!\not=\!\emptyset$, which proves the desired property for $p\! =\! 1$. Then we argue inductively on $p$. So assume that the property is proved for $p$. Note that 
$S_p\!\subseteq\!\overline{S_{p+1}}^{T_{\eta_p+1}}$, and 
$S_p\cap\bigcap_{1\leq i<p}~\overline{S_i}^{T_{\eta_i}}\cap D\cap O\!\not=\!\emptyset$, by induction assumption. Thus 
$$\overline{S_{p+1}}^{T_{\eta_p+1}}\cap\bigcap_{1\leq i\leq p}~\overline{S_i}^{T_{\eta_i}}\cap D\cap O
\!\not=\!\emptyset .$$ 
As $\bigcap_{1\leq i\leq p}~\overline{S_i}^{T_{\eta_i}}\cap D\cap O$ is $T_{\eta_p+1}$-open, 
$S_{p+1}\cap\bigcap_{1\leq i\leq p}~\overline{S_i}^{T_{\eta_i}}\cap D\cap O\!\not=\!\emptyset$.

\vfill\eject

\noindent (c) We use the notation before Lemma \ref{enum}. We enumerate $\{\xi_i\mid\xi_i\!\geq\! 1\}$ in an increasing way by $\{\eta_i\mid 1\!\leq\! i\!\leq\! p\}$, which means that we forget $\xi_1$ if it is $0$. As $\eta\!\geq\! 1$, $p\!\geq\! 1$. Note that we may assume that $s^{\eta_i+1}\!\subsetneqq\! s^{\eta_i}$ if $1\!\leq\! i\! <\! p$, by Lemma \ref{enum}. We set $S_i\! :=\! S_{s^{\eta_i}}$, for $1\!\leq\! i\!\leq\! p$. Note that $S_i\!\subseteq\!\overline{S_{i+1}}^{T_{\eta_i+1}}$ if $1\!\leq\! i\! <\! p$ since 
$s^{\eta_{i+1}}\ R^{\eta_i+1}\ s^{\eta_i}$. Thus 
$S_{s^\eta}\cap\bigcap_{1\leq\xi_i<\eta}~\overline{S_{s^{\xi_i}}}^{T_{\xi_i}}\cap O$ and 
$E\cap\bigcap_{1\leq\xi_i\leq\eta}~\overline{S_{s^{\xi_i}}}^{T_{\xi_i}}\cap O$ are $T_1$-dense in 
$\overline{S_{s^1}}^{T_1}\cap O$, by (b) and since $s^{\eta_1}\! =\! s^1$. But if $1\!\leq\!\rho\!\leq\!\eta$, then there is $1\!\leq\! i\!\leq\! n$ with $s^\rho\! =\! s^{\xi_i}$. And $\rho\!\leq\!\xi_i$ since 
$s^{\xi_i+1}\!\subsetneqq\! s^{\xi_i}$ if $1\!\leq\! i\! <\! n$. Thus we are done since 
$S_{s^{\eta}}\cap\bigcap_{1\leq\rho <\eta}~\overline{S_{s^\rho}}^{T_\rho}\! =\! 
S_{s^\eta}\cap\bigcap_{1\leq\xi_i <\eta}~\overline{S_{s^{\xi_i}}}^{T_{\xi_i}}$ and 
$\bigcap_{1\leq\rho\leq\eta}~\overline{S_{s^\rho}}^{T_\rho}\! =\!
\bigcap_{1\leq\xi_i\leq\eta}~\overline{S_{s^{\xi_i}}}^{T_{\xi_i}}$.\hfill{$\square$}

\subsection{$\!\!\!\!\!\!$ Some general effective facts}
 
\begin{lemm} \label{bisep} Let $1\!\leq\!\eta ,\xi\! <\!\omega_1^{\mbox{CK}}$, $X, Y$ be recursively presented Polish spaces, $A\!\in\!\Ana (X)\cap\borae$, $B\!\in\!\Ana (Y)\cap\boraxi$, and 
$C\!\in\!\Ana (X\!\times\! Y)$ disjoint from $A\!\times\! B$. Then there are $A'\!\in\!\Borel\cap\borae$, $B'\!\in\!\Borel\cap\boraxi$ such that $A'\!\times\! B'$ separates $A\!\times\! B$ from $C$. This also holds for the multiplicative classes.\end{lemm}
 
\noindent\bf Proof.\rm\ We argue as in the proof of Lemma 2.2 in [L-Z].\hfill{$\square$}

\begin{them} \label{effsep} Let $1\!\leq\!\eta ,\xi\! <\!\omega_1^{\mbox{CK}}$, $X, Y$ be recursively presented Polish spaces, and $A,B$ be disjoint $\Ana$ subsets of $X\!\times\! Y$. We assume that $A$ is separable from $B$ by a $(\borae\!\times\!\boraxi )_\sigma$ set. Then $A$ is separable from $B$ by a $\Borel\cap\big( (\Borel\cap\borae )\!\times\! (\Borel\cap\boraxi )\big)_\sigma$ set.\end{them}

\noindent\bf Proof.\rm\ We argue as in the proof of Theorem 2.3 in [L-Z].\hfill{$\square$}\medskip

 The next result is similar to Theorem 2.5 in [L-Z].

\begin{them} \label{kernel} Let $1\!\leq\!\eta ,\xi\! <\!\omega_1^{\mbox{CK}}$, $X, Y$ be recursively presented Polish spaces, and $A,B$ be disjoint $\Ana$ subsets of $X\!\times\! Y$. The following are equivalent:\smallskip

(a) the set $A$ cannot be separated from $B$ by a $(\borae\!\times\!\boraxi )_\sigma$ set.\smallskip

(b) the set $A$ cannot be separated from $B$ by a $\Borel\cap (\borae\!\times\!\boraxi )_\sigma$ set.\smallskip

(c) the set $A$ cannot be separated from $B$ by a $\boraone (T_\eta\!\times\! T_\xi )$ set.\smallskip

(d) $A\cap\overline{B}^{T_\eta\times T_\xi}\!\not=\!\emptyset$.\end{them}

\noindent\bf Proof.\rm\ Theorem \ref{effsep} implies that (a) is indeed equivalent to (b), and actually to the fact that $A$ cannot be separated from $B$ by a 
$\Borel\cap\big( (\Borel\cap\borae )\!\times\! (\Borel\cap\boraxi )\big)_\sigma$ set. By Theorem 1.A in [Lo], a $\Borel\cap\boraxi$ set is a countable union of $\Borel\cap\bormlxi$ sets, and thus 
$T_\xi$-open, if $\xi\!\geq\! 2$. Therefore (c) implies (a), and the converse is clear. It is also clear that (c) and (d) are equivalent.\hfill{$\square$}\medskip

 The following result is Lemma 3.3 in [Za], and is a consequence of Theorem \ref{kernel}.
 
\begin{them} \label{kernel2} Let $1\!\leq\!\xi ,\eta\! <\!\omega_1^{\mbox{CK}}$, $X, Y$ be recursively presented Polish spaces, and $A,B$ be disjoint $\Ana$ subsets of $X\!\times\! Y$. The following are equivalent:\smallskip

(a) The set $A$ cannot be separated from $B$ by a $\borme\!\times\!\bormxi$ set.\smallskip

(b) $B\cap (\overline{\mbox{proj}_X[A]}^{T_\eta}\!\times\!\overline{\mbox{proj}_Y[A]}^{T_\xi})\!\not=\!\emptyset$.
\end{them}

\vfill\eject

\section{$\!\!\!\!\!\!$ Countable unions of $\boraone\!\times\!\boraxi$ sets}\indent

 Let $Q\!\in\!\bormxi (2^\omega )\!\setminus\!\boraxi$. Then 
$P\! :=\! h[Q]\!\in\!\bormxi ([\subseteq ])\!\setminus\!\boraxi$ since $h$ is a homeomorphism.\medskip

\noindent\bf (A) The successor case\rm\medskip

 Assume that $\xi\! =\!\eta\!+\! 1$ is a countable ordinal. Theorem \ref{rep} gives a resolution family $(R^\rho )_{\rho\leq\eta}$. We set $\mathbb{X}\! :=\! [R^\eta ]$, 
$\mathbb{Y}\! :=\! [\subseteq ]$, 
$\mathbb{A}\! :=\!\{ (\beta ,\alpha )\!\in\!\mathbb{X}\!\times\!\mathbb{Y}\mid
\Pi (\beta )\! =\!\alpha\!\in\! P\}$ and 
$$\mathbb{B}\! :=\!\{ (\beta ,\alpha )\!\in\!\mathbb{X}\!\times\!\mathbb{Y}\mid
\Pi (\beta )\! =\!\alpha\!\notin\! P\} .$$ 
Note that $\mathbb{X}$ and $\mathbb{Y}$ are zero-dimensional Polish spaces, $\mathbb{A}$ is a closed subset of $\mathbb{X}\!\times\!\mathbb{Y}$, and $\mathbb{B}$ is a difference of two closed subsets of $\mathbb{X}\!\times\!\mathbb{Y}$, and disjoint from $\mathbb{A}$.

\begin{lem} \label{nonsep1x} The set $\mathbb{A}$ is not separable from $\mathbb{B}$ by a 
$(\boraone\!\times\!\boraxi )_\sigma$ subset of $\mathbb{X}\!\times\!\mathbb{Y}$.\end{lem}

\noindent\bf Proof.\rm\ We argue by contradiction, which gives a sequence $(O_n)_{n\in\omega}$ of $\boraone$ subsets of $[R^\eta ]$ and a sequence $(S_n)_{n\in\omega}$ of $\boraxi$ subsets of $[\subseteq ]$ such that 
$\mathbb{A}\!\subseteq\!\bigcup_{n\in\omega}~O_n\!\times\! S_n\!\subseteq\!\neg\mathbb{B}$. This implies that $P\! =\!\bigcup_{n\in\omega}~\Pi [O_n]\cap S_n$. As $\Pi^{-1}$ is 
$\boraxi$-measurable, $\Pi [O_n]\!\in\!\boraxi ([\subseteq ])$ and $P\!\in\!\boraxi ([\subseteq ])$, which is absurd.\hfill{$\square$}\medskip

\noindent\bf Proof of Theorem \ref{countunionrectsigma}.\rm\ The exactly part comes from Lemma \ref{nonsep1x}. Assume that (a) does not hold. In order to simplify the notation, we will asume that 
$\xi\! <\!\omega_1^{\mbox{CK}}$, $X$ and $Y$ are recursively presented and $A,B$ are $\Ana$, so that $N\! :=\! A\cap\overline{B}^{T_1\times T_\xi}$ is a nonempty (by Theorem \ref{kernel}) 
$\Ana$ (as in the proof of Lemma \ref{top}.(a)) subset of $X\!\times\! Y$.\medskip

 We set ${\cal I}\! :=\!\{ s\!\in\! 2^{<\omega}\mid N^{R^\eta}_s\cap\Pi^{-1}(P)\!\not=\!\emptyset\}$. As $B$ is not empty, we may assume that $P\!\not=\!\emptyset$. In particular, $\emptyset\!\in\! {\cal I}$. We construct, for $s\!\in\! 2^{<\omega}$,\medskip

\noindent - $x_s\!\in\! X$ and $X_s\!\in\!\Boraone (X)$,\smallskip

\noindent - $y_s\!\in\! Y$ and $Y_s\!\in\!\Boraone (Y)$,\smallskip

\noindent - $S_s\!\in\!\Ana (X\!\times\! Y)$.\medskip

 We want these objects to satisfy the following conditions:
$$\begin{array}{ll}
& (1)~\left\{\!\!\!\!\!\!\!\!
\begin{array}{ll}
& \overline{X_t}\!\subseteq\! X_s\mbox{ if }s~R^\eta ~t\wedge s\!\not=\! t\cr 
& \overline{Y_t}\!\subseteq\! Y_s\mbox{ if }s~R^0~t\wedge s\!\not=\! t\cr
& S_t\!\subseteq\! S_s\mbox{ if }s~R^\eta ~t\wedge (s,t\!\in\! {\cal I}\vee s,t\!\notin\! {\cal I})
\end{array}
\right.\cr\cr
& (2)~x_s\!\in\! X_s\wedge y_s\!\in\! Y_s\wedge 
(x_s,y_s)\!\in\! S_s\!\subseteq\! (X_s\!\times\! Y_s)\cap\Omega_{X\times Y}\cr\cr
& (3)~\mbox{diam}(X_s),\mbox{diam}(Y_s),\mbox{diam}_{\mbox{GH}}(S_s)\!\leq\! 2^{-\vert s\vert}\cr\cr
& (4)~S_s\!\subseteq\!\left\{\!\!\!\!\!\!\!\!
\begin{array}{ll}
& N\mbox{ if }s\!\in\! {\cal I}\cr
& B\mbox{ if }s\!\notin\! {\cal I}
\end{array}
\right.\cr\cr
& (5)~\mbox{proj}_Y[S_t]\!\subseteq\!\overline{\mbox{proj}_Y[S_s]}^{T_\rho}\mbox{ if }s~R^\rho ~t\wedge 
1\!\leq\!\rho\!\leq\!\eta
\end{array}$$

 Assume that this is done. Let $\beta\!\in\!\mathbb{X}$. Note that $\beta (k)~R^\eta ~\beta (k\! +\! 1)$ for each $k\!\in\!\omega$. By (1), 
$$\overline{X_{\beta (k+1)}}\!\subseteq\! X_{\beta (k)}.$$ 
Thus $(\overline{X_{\beta (k)}})_{k\in\omega}$ is a decreasing sequence of nonempty closed subsets of $X$ with vanishing diameters. We define 
${\{ f(\beta )\}\! :=\!\bigcap_{k\in\omega}~\overline{X_{\beta (k)}}\! =\!
\bigcap_{k\in\omega}~X_{\beta (k)}}$, so that 
$f(\beta )\! =\!\mbox{lim}_{k\rightarrow\infty}~x_{\beta (k)}$ and $f$ is continuous.\medskip

 Now let $\alpha\!\in\!\mathbb{Y}$. By (1), $\overline{Y_{\alpha (k+1)}}\!\subseteq\! Y_{\alpha (k)}$. Thus $(\overline{Y_{\alpha (k)}})_{k\in\omega}$ is a decreasing sequence of nonempty closed subsets of $Y$ with vanishing diameters. We define 
${\{ g(\alpha )\}\! :=\!\bigcap_{k\in\omega}~\overline{Y_{\alpha (k)}}\! =\!
\bigcap_{k\in\omega}~Y_{\alpha (k)}}$, so that 
$g(\alpha )\! =\!\mbox{lim}_{k\rightarrow\infty}~y_{\alpha (k)}$ and 
$g\! :\!\mathbb{Y}\!\rightarrow\! Y$ is continuous.\medskip

 Let $(\beta ,\alpha )\!\in\!\mathbb{A}$. Note that $\beta (k)\!\in\! {\cal I}$ for each $k\!\in\!\omega$. By (1)-(4), $(S_{\beta (k)})_{k\in\omega}$ is a decreasing sequence of nonempty clopen subsets of $N\cap\Omega_{X\times Y}$ with vanishing GH-diameters. We set 
${\{ F(\beta )\}\! :=\!\bigcap_{k\in\omega}~S_{\beta (k)}}$. Note that $(x_{\beta (k)},y_{\beta (k)})$ converge to $F(\beta )$ for ${\it\Sigma}_{X^2}$, and thus ${\it\Sigma}_X^2$. So their limit is 
$\big( f(\beta ),g(\alpha )\big)$, which is therefore in $N\!\subseteq\! A$, showing that 
$\mathbb{A}\!\subseteq\! (f\!\times\! g)^{-1}(A)$.\medskip

 Let $(\beta ,\alpha )\!\in\!\mathbb{B}$. As $\Pi^{-1}(P)$ is a closed subset of $[R^\eta ]$, there is $k_0\!\in\!\omega$ such that $\beta (k)\!\notin\! {\cal I}$ if $k\!\geq\! k_0$. By (1)-(4), 
$(S_{\beta (k)})_{k\geq k_0}$ is a decreasing sequence of nonempty clopen subsets of 
${B\cap\Omega_{X\times Y}}$ with vanishing GH-diameters, and we define 
$\{ G(\beta )\}\! :=\!\bigcap_{k\geq k_0}~S_{\beta (k)}$. Note that $(x_{\beta (k)},y_{\beta (k)})$ converge to $G(\beta )$. So their limit is $\big( f(\beta ),g(\alpha )\big)$, which is therefore in $B$, showing that $\mathbb{B}\!\subseteq\! (f\!\times\! g)^{-1}(B)$.\medskip

 Let us prove that the construction is possible. Let 
$(x_\emptyset ,y_\emptyset )\!\in\! N\cap\Omega_{X\!\times\! Y}$, and 
$X_\emptyset ,Y_\emptyset\!\in\!\Boraone$ with diameter at most $1$ such that 
$(x_\emptyset ,y_\emptyset )\!\in\! X_\emptyset\!\times\! Y_\emptyset$, as well as 
$S_\emptyset\!\in\!\Ana (X\!\times\! Y)$ with GH-diameter at most $1$ and 
$(x_\emptyset ,y_\emptyset )\!\in\! S_\emptyset\!\subseteq\! 
N\cap (X_\emptyset\!\times\! Y_\emptyset )\cap\Omega_{X\times Y}$. Assume that our objects satisfying (1)-(5) are constructed up to the length $l$, which is the case for $l\! =\! 0$. So let 
$s\!\in\! 2^{l+1}$.\medskip

\noindent\bf Claim\it\ The set $\mbox{proj}_Y[S_{s^\eta}]\cap
\bigcap_{1\leq\rho <\eta}~\overline{\mbox{proj}_Y[S_{s^\rho}]}^{T_\rho}\cap Y_{s^0}$ is $T_1$-dense in 
$\overline{\mbox{proj}_Y[S_{s^1}]}\cap Y_{s^0}$ if $\eta\!\geq\! 1$.\rm\medskip

 Indeed, we apply Lemma \ref{top}.(c) to $E\! :=\! Y$ and $O\! :=\! Y_{s^0}$.\hfill{$\diamond$}\medskip

 Note that $s^1\!\subseteq\! s^0\!\subsetneqq\! s$ and $s^1~R^1~s^0$, so that 
$\mbox{proj}_Y[S_{s^0}]\!\subseteq\!\overline{\mbox{proj}_Y[S_{s^1}]}$. Thus 
$y_{s^0}\!\in\!\overline{\mbox{proj}_Y[S_{s^1}]}\cap Y_{s^0}$. This shows that 
$I\! :=\!\mbox{proj}_Y[S_{s^\eta}]\cap
\bigcap_{1\leq\rho <\eta}~\overline{\mbox{proj}_Y[S_{s^\rho}]}^{T_\rho}\cap Y_{s^0}$ is not empty, even if $\eta\! =\! 0$.\medskip

\noindent\bf Case 1\rm\ $s\!\notin\! {\cal I}$\medskip

\noindent\bf 1.1\rm\ If $s^\eta\!\notin\! {\cal I}$, then we choose $y_s\!\in\! I$, $x_s\!\in\! X_{s^\eta}$ with $(x_s,y_s)\!\in\! S_{s^\eta}$, $X_s,Y_s\!\in\!\Boraone$ with diameter at most $2^{-l-1}$ such that $(x_s,y_s)\!\in\! X_s\!\times\! Y_s\!\subseteq\!\overline{X_s}\!\times\!\overline{Y_s}\!\subseteq\! 
X_{s^\eta}\!\times\! Y_{s^0}$, and also $S_s\!\in\!\Ana (X\!\times\! Y)$ with GH-diameter at most 
$2^{-l-1}$ such that $(x_s,y_s)\!\in\! S_s\!\subseteq\! S_{s^\eta}\cap
\big( X_s\!\times\! (\bigcap_{1\leq\rho <\eta}~\overline{\mbox{proj}_Y[S_{s^\rho}]}^{T_\rho}\cap Y_s)\big)$. If $s~R^\eta ~t$ and $s\!\not=\! t$, then $s~R^0~t^\eta ~R^\eta ~t$, so that $s~R^\eta ~t^\eta$, by Lemma \ref{extdist}. This implies that $\overline{X_t}\!\subseteq\! X_s$ and 
$\mbox{proj}_Y[S_{t^\eta}]\!\subseteq\!\overline{\mbox{proj}_Y[S_s]}^{T_\eta}$. Thus 
$\mbox{proj}_Y[S_t]\!\subseteq\!\overline{\mbox{proj}_Y[S_s]}^{T_\eta}$. If moreover $s\!\notin\! {\cal I}$, then 
$t^\eta\!\notin\! {\cal I}$ since $s~R^\eta ~t^\eta$. Thus $S_{t^\eta}\!\subseteq\! S_s$ and 
$S_t\!\subseteq\! S_s$. Similarly, $\overline{Y_t}\!\subseteq\! Y_s$ if $s~R^0~t$ and 
$s\!\not=\! t$ (this is simpler). If $1\!\leq\!\rho\! <\!\eta$, $s~R^\rho ~t$ and $s\!\not=\! t$, then 
$s~R^\rho ~t^\rho$, $\mbox{proj}_Y[S_{t^\rho}]\!\subseteq\!\overline{\mbox{proj}_Y[S_s]}^{T_\rho}$ and 
$\mbox{proj}_Y[S_t]\!\subseteq\!\overline{\mbox{proj}_Y[S_s]}^{T_\rho}$.\medskip

\noindent\bf 1.2\rm\ If $s^\eta\!\in\! {\cal I}$, then we choose $y\!\in\! I$, and $x\!\in\! X_{s^\eta}$ with 
$(x,y)\!\in\! S_{s^\eta}$. Note that 
$$(x,y)\!\in\!\overline{B}^{T_1\times T_\xi}\cap\big( X_{s^\eta}\!\times\! 
(\bigcap_{1\leq\rho\leq\eta}~\overline{\mbox{proj}_Y[S_{s^\rho}]}^{T_\rho}\cap Y_{s^0})\big) .$$ 

\vfill\eject

 This gives $(x_s,y_s)\!\in\! B\cap\big( X_{s^\eta}\!\times\! 
(\bigcap_{1\leq\rho\leq\eta}~\overline{\mbox{proj}_Y[S_{s^\rho}]}^{T_\rho}\cap Y_{s^0})\big)\cap
\Omega_{X\times Y}$. We choose $X_s,Y_s\!\in\!\Boraone$ with diameter at most $2^{-l-1}$ such that $(x_s,y_s)\!\in\! X_s\!\times\! Y_s\!\subseteq\!\overline{X_s}\!\times\!\overline{Y_s}\!\subseteq\! 
X_{s^\eta}\!\times\! Y_{s^0}$, and $S_s\!\in\!\Ana (X\!\times\! Y)$ with GH-diameter at most 
$2^{-l-1}$ such that $(x_s,y_s)\!\in\! S_s\!\subseteq\! B\cap
\big( X_s\!\times\! (\bigcap_{1\leq\rho\leq\eta}~\overline{\mbox{proj}_Y[S_{s^\rho}]}^{T_\rho}\cap Y_s)\big)\cap\Omega_{X\times Y}$. As above, we check that these objects are as required.\medskip

\noindent\bf Case 2\rm\ $s\!\in\! {\cal I}$\medskip

 Note that $s^\eta\!\in\! {\cal I}$. We argue as in 1.1.\hfill{$\square$}\medskip

\noindent\bf (B) The limit case\rm\medskip

 Assume that $\xi$ is an infinite limit ordinal. We indicate the differences with the successor case. Theorem \ref{replim} gives a uniform resolution family $(R^\rho )_{\rho\leq\xi}$. We set 
$\mathbb{X}\! :=\! [R^\xi ]$, $\mathbb{Y}\! :=\! [\subseteq ]$, 
$$\mathbb{A}\! :=\!\{ (\gamma ,\beta )\!\in\!\mathbb{X}\!\times\!\mathbb{Y}\mid
\Pi (\gamma )\! =\!\beta\!\in\! P\}\mbox{,}$$ 
and $\mathbb{B}\! :=\!\{ (\gamma ,\beta )\!\in\!\mathbb{X}\!\times\!\mathbb{Y}\mid
\Pi (\gamma )\! =\!\beta\!\notin\! P\}$.\medskip

\noindent\bf Proof of Theorem \ref{countunionrectsigma}.\rm\  This time, 
${\cal I}\! :=\!\{ s\!\in\! 2^{<\omega}\mid N^{R^\xi}_s\cap\Pi^{-1}(P)\!\not=\!\emptyset\}$. If 
$s\!\in\! 2^{<\omega}$, then we set, as in the proof of Theorem 2.4.4 in [L3], 
$\xi (s)\! :=\!\mbox{max}\{\xi_{h_{R^\xi}(t)+1}\mid t\!\subseteq\! s\}$. Note that $\xi (t)\!\leq\!\xi (s)$ if 
$t\!\subseteq\! s$. Conditions (1) and (5) become
$$\begin{array}{ll}
& (1')~\left\{\!\!\!\!\!\!\!\!
\begin{array}{ll}
& \overline{X_t}\!\subseteq\! X_s\mbox{ if }s~R^\xi ~t\wedge s\!\not=\! t\cr 
& \overline{Y_t}\!\subseteq\! Y_s\mbox{ if }s~R^0~t\wedge s\!\not=\! t\cr
& S_t\!\subseteq\! S_s\mbox{ if }s~R^\xi ~t\wedge (s,t\!\in\! {\cal I}\vee s,t\!\notin\! {\cal I})
\end{array}
\right.\cr\cr
& (5')~\mbox{proj}_Y[S_t]\!\subseteq\!\overline{\mbox{proj}_Y[S_s]}^{T_\rho}\mbox{ if }s~R^\rho ~t\wedge 
1\!\leq\!\rho\!\leq\!\xi (s)
\end{array}$$
The next claim and the remark after it were already present in the proof of Theorem 2.4.4 in [L3].\medskip

\noindent\bf Claim 1\rm ~Assume that $s^\rho\!\not=\! s^\xi$. Then $\rho\! +\! 1\!\leq\!\xi (s^{\rho +1})$.\medskip 

 We argue by contradiction. We get $\rho\! +\! 1\! >\!\rho\!\geq\!\xi (s^{\rho +1})\!\geq\!\xi_{h_{R^\xi}(s^{\xi})+1}\! =\!\xi_{h_{R^\xi}(s)}$. As $s^\rho\ R^\rho\ s$, 
$s^\rho\ R^\xi\ s$ and $s^\rho\! =\! s^\xi$, which is absurd.\hfill{$\diamond$}\medskip

 Note that $\xi_{n-1}\! <\!\xi_{n-1}\! +\! 1\!\leq\!\xi (s^{\xi_{n-1}+1})\!\leq\!\xi (s)$. Thus 
$s^{\xi (s)}\! =\! s^{\xi}$.\medskip

\noindent\bf Claim 2\rm ~The set $\mbox{proj}_Y[S_{s^\xi}]\cap\bigcap_{1\leq\rho <\xi (s)}~
\overline{\mbox{proj}_Y[S_{s^\rho}]}^{T_\rho}\cap Y_{s^0}$ is $T_1$-dense in 
$\overline{\mbox{proj}_Y[S_{s^1}]}^{T_1}\cap Y_{s^0}$.\medskip 

 We conclude as in the successor case, using the facts that $\xi_k\!\geq\! 1$ and $\xi (.)$ is increasing.\hfill{$\square$}

\section{$\!\!\!\!\!\!$ $\bormone\!\times\!\bormxi$ sets}\indent

 We consider $P$ as in Section 3.\medskip

\noindent\bf (A) The successor case\rm\medskip

 Assume that $\xi\! =\!\eta\!+\! 1$ is a countable ordinal. Theorem \ref{rep} gives a resolution family $(R^\rho )_{\rho\leq\eta}$. We set $\mathbb{X}\! :=\! [R^\eta ]\!\oplus\!\Pi^{-1}(\neg P)$, 
$\mathbb{Y}\! :=\! [\subseteq ]\!\oplus\!\Pi^{-1}(\neg P)$,  
$$\mathbb{A}\! :=\!\big\{\big( (0,\beta ),(1,\gamma )\big)\!\in\!\mathbb{X}\!\times\!\mathbb{Y}\mid
\beta\! =\!\gamma\big\}\cup
\big\{\big( (1,\gamma ),(0,\alpha )\big)\!\in\!\mathbb{X}\!\times\!\mathbb{Y}\mid
\Pi (\gamma )\! =\!\alpha\big\}$$ 
and $\mathbb{B}\! :=\!\big\{\big( (0,\beta ),(0,\alpha )\big)\!\in\!\mathbb{X}\!\times\!\mathbb{Y}\mid
\Pi (\beta )\! =\!\alpha\!\in\! P\big\}$. Note that $\mathbb{X}$ and $\mathbb{Y}$ are zero-dimensional Polish spaces, $\mathbb{A}$ is a closed subset of 
$\mathbb{X}\!\times\!\mathbb{Y}$, and $\mathbb{B}$ is a closed subset of 
$\mathbb{X}\!\times\!\mathbb{Y}$ disjoint from $\mathbb{A}$.

\begin{lem} \label{nonsep1xrect} The set $\mathbb{A}$ is not separable from $\mathbb{B}$ by a 
$\bormone\!\times\!\bormxi$ subset of $\mathbb{X}\!\times\!\mathbb{Y}$.\end{lem}

\noindent\bf Proof.\rm\ Let $C\!\in\!\bormone (\mathbb{X})$ and $S\!\in\!\bormxi (\mathbb{Y})$ with $\mathbb{A}\!\subseteq\! C\!\times\! S$. Note that 
$C\cap (\{0\}\!\times\! [R^\eta ])\! =\!\{0\}\!\times\! C'$ for some $C'\!\in\!\bormone ([R^\eta ])$. Similarly, $S\cap (\{0\}\!\times\! [\subseteq ])\! =\!\{0\}\!\times\! S'$ for some $S'\!\in\!\bormxi ([\subseteq ])$. Let $\alpha\!\in\! [\subseteq ]\!\setminus\! P$, and $\beta\! :=\!\gamma\! :=\!\Pi^{-1}(\alpha )$. Then 
$\big( (0,\beta ),(1,\gamma )\big)\!\in\!\mathbb{A}$, so that $\beta\!\in\! C'$ and $\alpha\!\in\!\Pi [C']$. Similarly, 
$\big( (1,\gamma ),(0,\alpha )\big)\!\in\!\mathbb{A}$, so that $\alpha\!\in\! S'$. This shows that $[\subseteq ]\!\setminus\! P\!\subseteq\!\Pi [C']\cap S'$. As 
$P\!\notin\!\boraxi ([\subseteq ])$, there is $\alpha\!\in\!\Pi [C']\cap S'\cap P$, and $\big( (0,\beta ),(0,\alpha )\big)\!\in\!\mathbb{B}\cap (C\!\times\! S)$ if 
$\beta\! :=\!\Pi^{-1}(\alpha )$.\hfill{$\square$}\medskip

\noindent\bf Proof of Theorem \ref{rectpi}.\rm\ The exactly part comes from Lemma \ref{nonsep1xrect}. Assume that (a) does not hold. In order to simplify the notation, we will assume that $\xi\! <\!\omega_1^{\mbox{CK}}$, $X$ and $Y$ are recursively presented and $A,B$ are 
$\Ana$, so that $N\! :=\! B\cap (\overline{\mbox{proj}_X[A]}\!\times\!\overline{\mbox{proj}_Y[A]}^{T_\xi})$ is a nonempty $\Ana$ subset of $X\!\times\! Y$, by Theorem \ref{kernel2}.\medskip

 We set ${\cal I}\! :=\!\{ s\!\in\! 2^{<\omega}\mid N^{R^\eta}_s\cap\Pi^{-1}(P)\!\not=\!\emptyset\}$. As $A$ is not empty, we may assume that $P\!\not=\!\emptyset$. In particular, $\emptyset\!\in\! {\cal I}$. We define, for $t\!\in\! 2^{<\omega}$, $t_c\!\in\! 2$ by $t_c\! :=\!\chi_{\neg {\cal I}}(t)$. We construct\medskip

\noindent - $x_{\varepsilon ,s}\!\in\! X$ and $X_{\varepsilon ,s}\!\in\!\Boraone (X)$ when 
$(\varepsilon ,s)\!\in\! (\{ 0\}\!\times\! 2^{<\omega})\cup\big(\{ 1\}\!\times\! (\neg {\cal I})\big)$,\smallskip

\noindent - $y_{\varepsilon ,s}\!\in\! Y$ and $Y_{\varepsilon ,s}\!\in\!\Boraone (Y)$ when 
$(\varepsilon ,s)\!\in\! (\{ 0\}\!\times\! 2^{<\omega})\cup\big(\{ 1\}\!\times\! (\neg {\cal I})\big)$,\smallskip

\noindent - $S_{\varepsilon ,\varepsilon',s}\!\in\!\Ana (X\!\times\! Y)$ when 
$(\varepsilon ,\varepsilon',s)\!\in\! 2^2\!\times\! 2^{<\omega}$, ($\varepsilon\!\not=\!\varepsilon'\wedge 
s\!\notin\! {\cal I}$) or $(\varepsilon\! =\!\varepsilon'\! =\! 0\wedge s\!\in\! {\cal I})$.\medskip

 We want these objects to satisfy the following conditions:
$$\begin{array}{ll}
& (1)~\left\{\!\!\!\!\!\!\!\!
\begin{array}{ll}
& \overline{X_{\varepsilon ,t}}\!\subseteq\! X_{\varepsilon ,s}\mbox{ if }s~R^\eta ~t\wedge 
s\!\not=\! t\cr 
& \overline{Y_{0,t}}\!\subseteq\! Y_{0,s}\mbox{ if }s~R^0~t\wedge s\!\not=\! t\cr
& \overline{Y_{1,t}}\!\subseteq\! Y_{1,s}\mbox{ if }s~R^\eta ~t\wedge s\!\not=\! t\cr 
& S_{\varepsilon ,\varepsilon',t}\!\subseteq\! S_{\varepsilon ,\varepsilon',s}\mbox{ if }s~R^\eta ~t
\end{array}
\right.\cr\cr
& (2)~x_{\varepsilon ,s}\!\in\! X_{\varepsilon ,s}\wedge y_{\varepsilon ,s}\!\in\! Y_{\varepsilon ,s}\wedge (x_{\varepsilon ,s},y_{\varepsilon',s})\!\in\! S_{\varepsilon ,\varepsilon',s}\!\subseteq\! 
(X_{\varepsilon ,s}\!\times\! Y_{\varepsilon',s})\cap\Omega_{X\times Y}\cr\cr
& (3)~\mbox{diam}(X_{\varepsilon ,s}),\mbox{diam}(Y_{\varepsilon ,s}),
\mbox{diam}_{\mbox{GH}}(S_{\varepsilon ,\varepsilon',s})\!\leq\! 2^{-\vert s\vert}\cr\cr
& (4)~S_{\varepsilon ,\varepsilon',s}\!\subseteq\!\left\{\!\!\!\!\!\!\!\!
\begin{array}{ll}
& N\mbox{ if }s\!\in\! {\cal I}\cr
& A\mbox{ if }s\!\notin\! {\cal I}
\end{array}
\right.\cr\cr
& (5)~\mbox{proj}_Y[S_{t_c,0,t}]\!\subseteq\!\overline{\mbox{proj}_Y[S_{s_c,0,s}]}^{T_\rho}\mbox{ if }s~R^\rho ~t\wedge 1\!\leq\!\rho\!\leq\!\eta
\end{array}$$
Assume that this is done. Let $(0,\gamma )\!\in\!\mathbb{X}$. Note that 
$\gamma (k)~R^\eta ~\gamma (k\! +\! 1)$ for each $k\!\in\!\omega$. By (1), 
$\overline{X_{0,\gamma (k+1)}}\!\subseteq\! X_{0,\gamma (k)}$. Thus 
$(\overline{X_{0,\gamma (k)}})_{k\in\omega}$ is a decreasing sequence of nonempty closed subsets of $X$ with vanishing diameters. We define 
${\{ f(0,\gamma )\}\! :=\!\bigcap_{k\in\omega}~\overline{X_{0,\gamma (k)}}\! =\!
\bigcap_{k\in\omega}~X_{0,\gamma (k)}}$, so that 
$$f(0,\gamma )\! =\!\mbox{lim}_{k\rightarrow\infty}~x_{0,\gamma (k)}$$ 
and $f$ is continuous on $\{ 0\}\!\times\! [R^\eta ]$. Now let $(1,\gamma )\!\in\!\mathbb{X}$. Note that moreover that there is $k_\gamma\!\in\!\omega$ minimal such that 
$\gamma (k)\!\notin\! {\cal I}$ if $k\!\geq\! k_\gamma$. We define $f(1,\gamma )$ similarly, using 
$(\overline{X_{1,\gamma (k)}})_{k\geq k_\gamma}$. Note that $f$ is continuous on 
$\{ 1\}\!\times\!\Pi^{-1}(\neg\mathbb{P})$ since $k_{\gamma'}\! =\! k_\gamma$ if 
$\gamma'\!\in\! N^{R^\eta}_{\gamma (k_\gamma )}$.\medskip

 Now let $(0,\alpha )\!\in\!\mathbb{Y}$. By (1), 
$\overline{Y_{0,\alpha (k+1)}}\!\subseteq\! Y_{0,\alpha (k)}$. Thus 
$(\overline{Y_{0,\alpha (k)}})_{k\in\omega}$ is a decreasing sequence of nonempty closed subsets of $Y$ with vanishing diameters. We define 
$${\{ g(0,\alpha )\}\! :=\!\bigcap_{k\in\omega}~\overline{Y_{0,\alpha (k)}}\! =\!
\bigcap_{k\in\omega}~Y_{0,\alpha (k)}}\mbox{,}$$ 
so that $g(0,\alpha )\! =\!\mbox{lim}_{k\rightarrow\infty}~y_{0,\alpha (k)}$. We define 
$g(1,\gamma )$ like $f(1,\gamma )$, so that $g\! :\!\mathbb{Y}\!\rightarrow\! Y$ is continuous.\medskip

 Assume that $\big( (0,\gamma ),(1,\gamma )\big)\!\in\!\mathbb{A}$. As $\Pi^{-1}(P)$ is a closed subset of $[R^\eta ]$, there is $k_0\!\in\!\omega$ such that $\gamma (k)\!\notin\! {\cal I}$ if 
$k\!\geq\! k_0$. By (1)-(4), $(S_{0,1,\gamma (k)})_{k\geq k_0}$ is a decreasing sequence of nonempty clopen subsets of $A\cap\Omega_{X\times Y}$ with vanishing GH-diameters. We set 
${\{ F(\gamma )\}\! :=\!\bigcap_{k\geq k_0}~S_{0,1,\gamma (k)}}$. Note that 
$(x_{0,\gamma (k)},y_{1,\gamma (k)})$ converge to $F(\gamma )$ for ${\it\Sigma}_{X^2}$, and thus ${\it\Sigma}_X^2$. So their limit is $\big( f(0,\gamma ),g(1,\gamma )\big)$, which is therefore in $A$. If now $\big( (1,\gamma ),(0,\alpha )\big)\!\in\!\mathbb{A}$, then we argue similarly, showing that $\mathbb{A}\!\subseteq\! (f\!\times\! g)^{-1}(A)$.\medskip

 Let $\big( (0,\gamma ),(0,\alpha )\big)\!\in\!\mathbb{B}$. Note that $\gamma (k)\!\in\! {\cal I}$ for each $k\!\in\!\omega$. By (1)-(4), $(S_{0,0,\gamma (k)})_{k\in\omega}$ is a decreasing sequence of nonempty clopen subsets of ${N\cap\Omega_{X\times Y}}$ with vanishing GH-diameters, and we define $\{ G(\gamma )\}\! :=\!\bigcap_{k\in\omega}~S_{0,0,\gamma (k)}$. Note that 
$(x_{0,\gamma (k)},y_{0,\gamma (k)})$ converge to $G(\gamma )$. So their limit is 
$\big( f(0,\gamma ),g(0,\alpha )\big)$, which is therefore in $N\!\subseteq\! B$, showing that 
$\mathbb{B}\!\subseteq\! (f\!\times\! g)^{-1}(B)$.\medskip

 Let us prove that the construction is possible. Let 
$(x_{0,\emptyset},y_{0,\emptyset})\!\in\! N\cap\Omega_{X\!\times\! Y}$, and 
$X_{0,\emptyset} ,Y_{0,\emptyset}\!\in\!\Boraone$ with diameter at most $1$ such that 
$(x_{0,\emptyset},y_{0,\emptyset})\!\in\! X_{0,\emptyset}\!\times\! Y_{0,\emptyset}$, as well as 
$S_{0,0,\emptyset }\!\in\!\Ana (X\!\times\! Y)$ with GH-diameter at most $1$ and 
$(x_{0,\emptyset},y_{0,\emptyset})\!\in\! S_{0,0,\emptyset}\!\subseteq\! 
N\cap (X_{0,\emptyset}\!\times\! Y_{0,\emptyset})\cap\Omega_{X\times Y}$. Assume that our objects satisfying (1)-(5) are constructed up to the length $l$, which is the case for $l\! =\! 0$. So let 
$s\!\in\! 2^{l+1}$.\medskip

\noindent\bf Claim\it\ The set $\mbox{proj}_Y[S_{s^\eta_c,0,s^\eta}]\cap\bigcap_{1\leq\rho <\eta}~
\overline{\mbox{proj}_Y[S_{s^\rho_c,0,s^\rho}]}^{T_\rho}\cap Y_{0,s^0}$ is $T_1$-dense in 
$\overline{\mbox{proj}_Y[S_{s^1_c,0,s^1}]}\cap Y_{0,s^0}$ if $\eta\!\geq\! 1$.\rm\medskip

 As in the proof of Theorem \ref{countunionrectsigma}, we infer that 
$$I\! :=\!\mbox{proj}_Y[S_{s^\eta_c,0,s^\eta}]\cap\bigcap_{1\leq\rho <\eta}~
\overline{\mbox{proj}_Y[S_{s^\rho_c,0,s^\rho}]}^{T_\rho}\cap Y_{0,s^0}$$ 
is not empty.\medskip

\noindent\bf Case 1\rm\ $s\!\notin\! {\cal I}$\medskip

\noindent\bf 1.1\rm\ $s^\eta\!\notin\! {\cal I}$\medskip

 Note that $s^\eta_c\! =\! 1$. We choose $y_{0,s}\!\in\! I$, $x_{1,s}\!\in\! X_{1,s^\eta}$ with 
$(x_{1,s},y_{0,s})\!\in\! S_{1,0,s^\eta}$, $X_{1,s},Y_{0,s}\!\in\!\Boraone$ with diameter at most 
$2^{-l-1}$ such that $(x_{1,s},y_{0,s})\!\in\! X_{1,s}\!\times\! Y_{0,s}\!\subseteq\!
\overline{X_{1,s}}\!\times\!\overline{Y_{0,s}}\!\subseteq\! X_{1,s^\eta}\!\times\! Y_{0,s^0}$, 
and also $S_{1,0,s}\!\in\!\Ana (X\!\times\! Y)$ with GH-diameter at most $2^{-l-1}$ such that 
$$(x_{1,s},y_{0,s})\!\in\! S_{1,0,s}\!\subseteq\! S_{1,0,s^\eta}\cap\big( X_{1,s}\!\times\! 
(\bigcap_{1\leq\rho <\eta}~\overline{\mbox{proj}_Y[S_{s^\rho_c,0,s^\rho}]}^{T_\rho}\cap Y_{0,s})\big) .$$
As in the proof of Theorem \ref{countunionrectsigma}, we check that these objects are as required. We also set 
$$(x_{0,s},y_{1,s})\! :=\! (x_{0,s^\eta},y_{1,s^\eta})\mbox{,}$$ 
choose $X_{0,s},Y_{1,s}\!\in\!\Boraone$ with diameter at most $2^{-l-1}$ such that 
$$(x_{0,s},y_{1,s})\!\in\! X_{0,s}\!\times\! Y_{1,s}\!\subseteq\!
\overline{X_{0,s}}\!\times\!\overline{Y_{1,s}}\!\subseteq\! X_{0,s^\eta}\!\times\! Y_{1,s^\eta}\mbox{,}$$ 
and also $S_{0,1,s}\!\in\!\Ana (X\!\times\! Y)$ with GH-diameter at most $2^{-l-1}$ such that 
$$(x_{0,s},y_{1,s})\!\in\! S_{0,1,s}\!\subseteq\! S_{0,1,s^\eta}\cap (X_{0,s}\!\times\! Y_{1,s}).$$

\vfill\eject

\noindent\bf 1.2\rm\ $s^\eta\!\in\! {\cal I}$\medskip

 We choose $y\!\in\! I$, and $x\!\in\! X$ with $(x,y)\!\in\! S_{0,0,s^\eta}$. Note that 
$$y\!\in\!\overline{\mbox{proj}_Y[A]}^{T_\xi}\cap
\bigcap_{1\leq\rho\leq\eta}~\overline{\mbox{proj}_Y[S_{s^\rho_c,0,s^\rho}]}^{T_\rho}\cap Y_{0,s^0}.$$
This gives $y'\!\in\!\mbox{proj}_Y[A]\cap\bigcap_{1\leq\rho\leq\eta}~
\overline{\mbox{proj}_Y[S_{s^\rho_c,0,s^\rho}]}^{T_\rho}\cap Y_{0,s^0}$, $x'\!\in\! X$ with 
$$(x',y')\!\in\! A\cap\big( X\!\times\! (\bigcap_{1\leq\rho\leq\eta}~
\overline{\mbox{proj}_Y[S_{s^\rho_c,0,s^\rho}]}^{T_\rho}\cap Y_{0,s^0})\big)\mbox{,}$$ 
and also $(x_{1,s},y_{0,s})\!\in\! A\cap\big( X\!\times\! 
(\bigcap_{1\leq\rho\leq\eta}~\overline{\mbox{proj}_Y[S_{s^\rho_c,0,s^\rho}]}^{T_\rho}\cap Y_{0,s^0})\big)
\cap\Omega_{X\times Y}$. We choose $X_{1,s},Y_{0,s}$ in $\Boraone$ with diameter at most 
$2^{-l-1}$ such that $(x_{1,s},y_{0,s})\!\in\! X_{1,s}\!\times\! Y_{0,s}\!\subseteq\!
\overline{X_{1,s}}\!\times\!\overline{Y_{0,s}}\!\subseteq\! X\!\times\! Y_{0,s^0}$, and 
$S_{1,0,s}\!\in\!\Ana (X\!\times\! Y)$ with GH-diameter at most $2^{-l-1}$ such that 
$$(x_{1,s},y_{0,s})\!\in\! S_{1,0,s}\!\subseteq\! A\cap
\big( X_{1,s}\!\times\! (\bigcap_{1\leq\rho\leq\eta}~
\overline{\mbox{proj}_Y[S_{s^\rho_c,0,s^\rho}]}^{T_\rho}\cap Y_{0,s})\big)\cap\Omega_{X\times Y} .$$ 
As above, we check that these objects are as required.\medskip

 Note also that $(x_{0,s^\eta},y_{0,s^\eta})\!\in\! S_{0,0,s^\eta}$, so that 
$x_{0,s^\eta}\!\in\!\overline{\mbox{proj}_X[A]}\cap X_{0,s^\eta}$. This gives a point 
$x'$ of $\mbox{proj}_X[A]\cap X_{0,s^\eta}$, and $y'\!\in\! Y$ with $(x',y')\!\in\! A\cap (X_{0,s^\eta}\!\times\! Y)$, and $(x_{0,s},y_{1,s})\!\in\! A\cap (X_{0,s^\eta}\!\times\! Y)\cap\Omega_{X\times Y}$. We choose $X_{0,s},Y_{1,s}\!\in\!\Boraone$ with diameter at most $2^{-l-1}$ such that 
$$(x_{0,s},y_{1,s})\!\in\! X_{0,s}\!\times\! Y_{1,s}\!\subseteq\!
\overline{X_{0,s}}\!\times\!\overline{Y_{1,s}}\!\subseteq\! X_{0,s^\eta}\!\times\! Y\mbox{,}$$ 
and $S_{0,1,s}\!\in\!\Ana (X\!\times\! Y)$ with GH-diameter at most $2^{-l-1}$ such that 
$$(x_{0,s},y_{1,s})\!\in\! S_{0,1,s}\!\subseteq\! 
A\cap (X_{0,s}\!\times\! Y_{0,s})\cap\Omega_{X\times Y}.$$ 
As above, we check that these objects are as required.\medskip

\noindent\bf Case 2\rm\ $s\!\in\! {\cal I}$\medskip

 Note that $s^\eta\!\in\! {\cal I}$. We argue as in the first part of 1.1 to construct $x_{0,s}$, $y_{0,s}$, $X_{0,s}$, $Y_{0,s}$ and $S_{0,0,s}$.\hfill{$\square$}\medskip
 
\noindent\bf (B) The limit case\rm\medskip

 Assume that $\xi$ is an infinite limit ordinal. We indicate the differences with the successor case. Theorem \ref{replim} gives a uniform resolution family $(R^\rho )_{\rho\leq\xi}$. We set 
$\mathbb{X}\! :=\! [R^\xi ]\!\oplus\!\Pi^{-1}(\neg P)$, 
$$\mathbb{Y}\! :=\! [\subseteq ]\!\oplus\!\Pi^{-1}(\neg P)\mbox{,}$$  
$\mathbb{A}\! :=\!\big\{\big( (0,\beta ),(1,\gamma )\big)\!\in\!\mathbb{X}\!\times\!\mathbb{Y}\mid
\beta\! =\!\gamma\big\}\cup
\big\{\big( (1,\gamma ),(0,\alpha )\big)\!\in\!\mathbb{X}\!\times\!\mathbb{Y}\mid
\Pi (\gamma )\! =\!\alpha\big\}$ and 
$$\mathbb{B}\! :=\!\big\{\big( (0,\beta ),(0,\alpha )\big)\!\in\!\mathbb{X}\!\times\!\mathbb{Y}\mid
\Pi (\beta )\! =\!\alpha\!\in\! P\big\} .$$
\bf Proof of Theorem \ref{rectpi}.\rm\ Condition (1) becomes
$$(1')~\left\{\!\!\!\!\!\!\!\!
\begin{array}{ll}
& \overline{X_{\varepsilon ,t}}\!\subseteq\! X_{\varepsilon ,s}\mbox{ if }s~R^\xi ~t\wedge 
s\!\not=\! t\cr 
& \overline{Y_{0,t}}\!\subseteq\! Y_{0,s}\mbox{ if }s~R^0~t\wedge s\!\not=\! t\cr
& \overline{Y_{1,t}}\!\subseteq\! Y_{1,s}\mbox{ if }s~R^\xi ~t\wedge s\!\not=\! t\cr 
& S_{\varepsilon ,\varepsilon',t}\!\subseteq\! S_{\varepsilon ,\varepsilon',s}\mbox{ if }s~R^\xi ~t
\end{array}
\right.$$
\bf Claim 2\rm ~The set $\mbox{proj}_Y[S_{s^\xi_c,0,s^\xi}]\cap\bigcap_{1\leq\rho <\xi (s)}~
\overline{\mbox{proj}_Y[S_{s^\rho_c,0,s^\rho}]}^{T_\rho}\cap Y_{0,s^0}$ is $T_1$-dense in 
$$\overline{\mbox{proj}_Y[S_{s^1_c,0,s^1}]}^{T_1}\cap Y_{0,s^0}.$$
We conclude as in the successor case.\hfill{$\square$}  
 
\section{$\!\!\!\!\!\!$ References}

\noindent [D-SR]\ \ G. Debs and J. Saint Raymond, Borel liftings of Borel sets: 
some decidable and undecidable statements,~\it Mem. Amer. Math. Soc.\rm ~187, 876 (2007)

\noindent [K]\ \ A. S. Kechris,~\it Classical Descriptive Set Theory,~\rm 
Springer-Verlag, 1995

\noindent [K-S-T]\ \ A. S. Kechris, S. Solecki and S. Todor\v cevi\'c, Borel chromatic numbers,\ \it 
Adv. Math.\rm\ 141 (1999), 1-44

\noindent [Ku]\ \ C. Kuratowski, Sur une g\'en\'eralisation de la notion d'hom\'eomorphie,~\it Fund. Math.\rm ~22, 1 (1934), 206-220

\noindent [L1]\ \ D. Lecomte, On minimal non potentially closed subsets of the plane,\ \it Topology Appl.\rm\ 154, 1 (2007), 241-262

\noindent [L2]\ \ D. Lecomte, A dichotomy characterizing analytic graphs of uncountable Borel chromatic number in any dimension,~\it Trans. Amer. Math. Soc.\rm~361 (2009), 4181-4193

\noindent [L3]\ \ D. Lecomte, How can we recognize potentially $\bormxi$ subsets of the plane?,~\it  J. Math. Log.\ \rm  9, 1 (2009), 39-62

\noindent [L4]\ \ D. Lecomte, Potential Wadge classes,~\it\ Mem. Amer. Math. Soc.,\rm ~221, 1038 (2013)

\noindent [L-Z]\ \ D. Lecomte and M. Zelen\'y, Baire-class $\xi$ colorings: the first three levels,\ \it Trans. Amer. Math. Soc.\rm\ 366, 5 (2014), 2345-2373

\noindent [Lo]\ \ A. Louveau, A separation theorem for $\Ana$ sets,\ \it Trans. 
Amer. Math. Soc.\ \rm 260 (1980), 363-378

\noindent [Lo-SR]\ \ A. Louveau and J. Saint Raymond, Borel classes and closed games: 
Wadge-type and Hurewicz-type results,\ \it Trans. Amer. Math. Soc.\ \rm 304 (1987), 431-467

\noindent [M]\ \ B. D. Miller, The graph-theoretic approach to descriptive set theory,\ \it Bull. Symbolic Logic\ \rm 18, 4 (2012), 554-575

\noindent [Mo]\ \ Y. N. Moschovakis,~\it Descriptive set theory,~\rm North-Holland, 1980

\noindent [Za]\ \ R. Zamora, Separation of analytic sets by rectangles of low complexity,~\it Topology Appl.~\rm 182 (2015), 77-97
 
\end{document}